# Modelling Chaotic Data

by

Vincent A. Mellor

MSc. Dissertation submitted in partial fulfilment of Data Analysis, Networks and Non-linear Dynamics Course in Mathematics Department at the University of York, UK.

August 2002

# ABSTRACT

This paper extends the subjects dicussed in the Data Analysis and Dynamical Systems courses by looking at the subject of modelling data. This task is nontrivial as the underlying process could be non-linear. In the paper some common methods, including global and local polynomial fitting, are discussed in terms of their applicability, level of computation and accuracy. One example method, Measure based Reconstruction, has been investigated in greater detail and experimentation is carried out to evaluate the method.

In this project we shall be looking at the different ways one can model chaotic time series. The reason we are going to look at a range of methods is that different methods are "good" for different applications. As the "goodness" of a model is subjective to the task one wishes to do, we will investigate a selected models and compare the prediction to see how one goes about testing a model.

# ACKNOWLEDGEMENTS

During the project the following people have been there to listen, guide and sometimes check what I had done. These people include my project Supervisor, Dr Zaqueu Coelho, who was there to point me in the right direction and listen to my concerns about the project. Dr Jason Levesley for generally suggest more ways to investigate and explaining the concepts of certain methods when Zaq wasn't around. Ian Naden for generally taking my mind of the project, and working out numerous equations to prove it wasn't just me getting the method not to work. Of course most of the project could have been done without the "research assistants" at home who were there for formatting and other general sitting and listening of ideas for the project.

# TABLE OF CONTENTS

# 1 INTRODUCTION

The goal of modelling is to establish the equations of motion that describe the underlying process of the data in terms of meaningful quantities. These models try to reproduce the internal properties of the underlying system of a time series so that it can be identified and then possibly used to create future points. Ideally the model would be able to extract these properties from little more than the time series itself. However, with a great deal of methods either assumptions must be put on to the data or some other method for determining the variables of the fitting must be used.

When choosing a model many considerations must be taken. The type of model required will normally depend on its application. If rough estimates are needed quickly, a large order approximation of the underlying function will not be required and the best model might be a polynomial fitted to the trend of the data. However, if an accurate prediction was required, this method would be useless as the errors of fitting would be large. The choice of model can also be affected by the type of data that is going to be modelled. For example stochastic models are good at modelling Gaussian white noise, but in general they are not very good for real world data.

In this project the most common methods in many classes will be discussed to give an idea of what properties each class of models have. This will then be followed by an investigation of the Measure Based Reconstruction method, with some numerical experiments to test the method. These experiments will be a guide to those that should be carried out on methods as a test to see how they will behave with actual data. To do these experiments the programs listed at the end of this document will be used.

# 2 REVIEW OF THEORY

In this section some topics that have been covered in lectures are briefly outlined to aid the explanation of future topics. These topics will only be quickly discussed with enough detail to introduce them and give some understanding for future sections. Most of the definitions have been found in the lecture notes for the courses, although some have been adapted from the books [1] and [2]. If more information is required on any of the topics it is advised that you start by looking at these books.

## 2.1 DYNAMICAL SYSTEMS

A dynamical system is a system which evolves in time in a deterministic way. In discrete time, given a set of measurement

$$x_0, x_1, x_2, \ldots, x_N \quad (N \text{ very large})$$

one thinks that there exists a function $f$ which is generating the data, i.e.

$$x_1 = f(x_0), x_2 = f(x_1), \ldots, x_N = f(x_{N-1})$$

$$\text{in short } x_{n+1} = f(x_n), \ n = 0, 1, \ldots, N-1$$

In the continuous time case we have a "continuous" path $x_t$ for $t \in [0, T]$. In this case it is thought that $x_t$ comes from a solution to an ordinary differential equation, with initial condition $x_0$ and we write:

$$x_t = \varphi_t(x_0)$$

where $\varphi_t$ is a flow. Most measurements made in the real world happen in discrete time (even when the underlying model is continuous), for instance temperature, pressure, etc. for weather prediction.

A flow is the unique solution to the ordinary differential equation on which the continuous dynamical system is based with initial condition $x(0) = x_0$. The flow $\varphi_t$ will satisfy the following properties:
1. $\varphi_0 =$ Identity map ($\varphi_0(x) \equiv x$)
2. $\varphi_{t+s} = \varphi_t \circ \varphi_s$ (flow property).

This works due to the uniqueness of the solution for a given initial condition.

## 2.2 STATE SPACE

The state variables of a dynamical system are the variables that fully describe its behaviour. For example the motion of a body in space is described fully by its position and its velocity. The state space $X$ of a system is the set of all possible states of the system. It is assumed that $X$ is a subset of some (possibly multidimensional) Euclidean space $\mathbb{R}^m$. The number of necessary state variables is called the dimension of the system.

A state is a vector $(x_1, x_2, \ldots, x_m)$ of real numbers. Let **$x_1$** be some initial state. This is where observation of the system begins. A trajectory is a sequence of points

$$x_1, x_2 \ldots$$

satisfying the time evolution of the system or a solution to the differential equation with initial condition $x_1$.

## 2.3 CANTOR SETS

A set of real numbers is called connected if it is an interval. A set is called totally disconnected if it does not contain any non empty interval. A set $S$ is called perfect if it does not contain isolated points. An isolated point is a point such that there exists $\varepsilon = \varepsilon(x)$ such that $[x-\varepsilon, x+\varepsilon] \cap S = \{x\}$. A Cantor Set is a set of real numbers, totally disconnected and perfect.

## 2.4 METRIC SPACES

Let $M$ be a set. A distance or metric on $M$ is a real function $d: M \times M \to \mathbb{R}$ satisfying:
1. $0 \leq d(x, y) < \infty$,
2. $d(x, y) = 0 \Leftrightarrow x = y$,
3. $d(x, y) = d(y, x)$,
4. $d(x, y) \leq d(x, z) + d(z, y)$ for all $x, y \in M$.

If a real function $d: M \times M \to \mathbb{R}$ satisfies 1, 3 and 4 then it is called pseudo metric on $M$. Property 4 is called the triangle inequality.

**Example**

$\mathbb{R}^n$ with Euclidean distance $d(x, y) = \sqrt{(x_1 - y_1)^2 + \ldots + (x_n - y_n)^2}$.

A set $M$ with distance $d$ is called a metric space and shall be denoted by $(M, d)$. Given a point $x \in M$ and a sequence of points $x_n \in M$ we say that $x_n$ converges to $x$ if $d(x_n, n) \to 0$ as $n \to \infty$. This means that for any $\varepsilon > 0$ there exists $n_0 = n_0(\varepsilon)$ such that for all $n \geq n_0$ we have $d(x_n, x) \leq \varepsilon$.

A sequence $x_n \in M$ is called a Cauchy sequence of for any $\varepsilon > 0$ there exists $n_0 = n_0(\varepsilon)$ such that for all $m, n \geq n_0$ we have $d(x_m, x_n) \leq \varepsilon$. **Note**: Every convergent sequence in $M$ is a Cauchy sequence. A metric space $(M, d)$ in which every Cauchy sequence is convergent is called a complete metric. The example earlier is a complete metric.

## 2.5 MEASURES

Let $X$ be a set. A $\sigma$-algebra (or $\sigma$-field) of subsets of $X$ is a collection $\boldsymbol{B}$ of subsets of $X$ satisfying the following properties:

1. $X \in \boldsymbol{B}$
2. If $B \in \boldsymbol{B}$ then $B^c \in \boldsymbol{B}$
3. If $B_n \in \boldsymbol{B}$ for $n > 1$ then $\bigcup_{n=1}^{\infty} B_n \in \boldsymbol{B}$

We call the pair $(X, \boldsymbol{B})$ a measurable space.

A finite measure on $(X, \boldsymbol{B})$ is a function $\mu : B \to \mathbb{R}_+$ satisfying $\mu(\emptyset) = 0$ and $\mu\left(\bigcup_{n=1}^{\infty} B_n\right) = \sum_{n=1}^{\infty} \mu(B_n)$ whenever $\{B_n\}$ is a collection of subsets which are pairwise disjoint. A finite measure space is a triple $(X, \boldsymbol{B}, \mu)$ where $\mu$ is a finite measure on the measurable space $(X, \boldsymbol{B})$. We say that $(X, \boldsymbol{B}, \mu)$ is a probability space if $\mu(X) = 1$. We then say that $\mu$ is a probability measure on the measurable space $(X, \boldsymbol{B})$

## 2.6 ERGODIC THEOREM

Let $X_i$ be the random outcome corresponding to the repetition of a random experiment. Put $X_i = 1$ if the experiment in the $i^{th}$ attempt is successful and $X_i = 0$ if it is a failure. Let the probability of success be denoted by $p$ ($0 < p < 1$). So we have $P(X_i = 1) = p$ and $P(X_i = 0) = 1 - p$. By repetition we mean that $X_i$ are independent random variables

$$P(X_1 = i_1, ..., X_n = i_n) = \prod_{j=1}^{n} P(X_j = i_j)$$

The law of large numbers

$$\frac{X_1 + X_2 + ... + X_n}{n} \xrightarrow[n \to \infty]{} E(X_1) \text{ "}P \text{ almost surely"}$$

Note that $E(X_1) = 1 \cdot P(X_1 = 1) + 0 \cdot P(X_1 = 0) = p$. So the average number of successes is asymptotically equal to the probability of success.

Now suppose the data $X_1, X_2, ..., X_n$ comes from measuring a dynamical system

$$X_i = \mathbf{1}_A(f^{i-1}(x))$$

$$\text{where } \mathbf{1}_A(x) = 1 \text{ if } x \in A$$

$$= 0 \text{ if } x \notin A$$

is the indicator function of $A$.

$$\frac{X_1 + X_2 + ... + X_n}{n} = \frac{1}{n} \sum_{j=0}^{n-1} \mathbf{1}_A(f^j(x))$$

is the average number of times the orbit of $x$ visits a set $A$. In the limit when $n \to \infty$ we should get a spatial mean of $X_1 = \mathbf{1}_A(x)$. This is expressed by the word ergodicity. So if the dynamical system is ergodic we should expect

$$\frac{1}{n} \sum_{j=0}^{n-1} \mathbf{1}_A(f^j(x)) \xrightarrow[n \to \infty]{} E(\mathbf{1}_A) = P(A) \text{ "}P \text{ almost surely"}.$$

## 2.7 INVARIANTS

Invariant qualities are those who's value is not affected under time evolution. These qualities will have the same value irrespective of the details of the measurement process and of the reconstruction of the state space. This is strictly true only for ideal, noise free and infinitely long time series, but a good algorithm applied to an approximately noise-free and sufficiently long data set should yield results which are robust against small changes in parameters of the algorithm. Examples of invariant qualities are the Lyapunov exponent and the correlation dimension.

### 2.7.1 Lyapunov Exponents

Given $f: I \to I$ and $x \in I$ the Lyapunov exponent of $x$ is given by the limit:

$$\lambda_x = \lim_{n \to \infty} \frac{1}{n} \log \left| (f^n)'(x) \right|$$

when the limit exists. Note that from the chain rule

$$(f^n)'(x) = \prod_{j=0}^{n-1} f'(f^j(x))$$

is the product of the derivative of *f* at the positive orbit of *x*. This is the growth rate of the derivative. Positive Lyapunov exponents corresponds to expansion along the positive orbit of *x*. Since the orbit is confined to an interval (bounded interval) then it implies chaotic behaviour.

## 2.8 FOURIER FREQUENCIES

Given $x_0, x_1, ..., x_{n-1}$. Frequencies that are integer multiples of $\frac{1}{n}$ are known as Fourier frequencies. They are distinguished by the fact that a sinusoid at that frequency will complete a whole number of cycles. Such frequencies are said to be harmonic with respect to the span of the data.

## 2.9 MARKOV PROCESSES

If the model starts in a given microstate *i*, it will move to state *j* with transition probability $P_{j \leftarrow i}$ that does not depend on the previous history of the model. If we assume the model to be under some fairly general conditions, processes after the passage of a transient would produce states with a unique steady-state probability distribution. This steady-state probability $\pi_j$ is an eigenvector, with eigenvalue one, of the transition matrix:

$$\pi_j = \sum_i P_{j \leftarrow i} \pi_i.$$

The steady-state probabilities are unique if the matrix $P_{j \leftarrow i}$ is regular, which means that for some integer *n* all elements of $(P_{j \leftarrow i})^n$ are positive and non-zero. Physically, this restriction implies that it is always possible to go from one state to any other state in a finite number of steps. Exceptions are matrices that are block diagonal, for example

$$\begin{bmatrix} P_{1 \leftarrow 1} & P_{1 \leftarrow 2} & 0 & 0 \\ P_{2 \leftarrow 1} & P_{2 \leftarrow 2} & 0 & 0 \\ 0 & 0 & P_{3 \leftarrow 3} & P_{3 \leftarrow 4} \\ 0 & 0 & P_{4 \leftarrow 3} & P_{4 \leftarrow 4} \end{bmatrix}$$

Since there is no way of going from states 1 or 2 to 3 or 4 the stationary probability distribution will depend on whether one started with one of the first two states or one of the last two.

## 2.10 Poincaré Sections

Let $S$ be a $n-1$ dimensional surface transverse to the trajectories of a dynamical system. Consider a point $x_0 \in S$ at time $t=0$. As the trajectory starting at $x_0$ evolves it will eventually return to $S$ at the point $x_1$ after some time $t_1$. Consider all initial points on $S$. Define the Poincaré return map $P : S \to S$ by

$$x_{k+1} = P(x_k).$$

The Poincaré section is formed by considering the intersection of $S$ with the attractor in its state space. The normal choice for $S$ is a hyperplane for simplicity and that this hyperplane is transverse to the flow. This means that is not tangential.

## 2.11 Taken's Embedding Theorem (1981)

Let $M$ be a compact manifold of dimension $m$. For generic differentiable map $f : M \to M$ and generic differentiable function $g : M \to \mathbb{R}$ the following map:

$$\Phi_{(f,g)}(x) = (g(x), g(f(x)), ..., g(f^{2m}(x)))$$

is an embedding of $M$ into $\mathbb{R}^{2m+1}$.

Generic means that given a differentiable map there exists another one arbitrarily near that one for which the theorem applies. Embedding means that $\Phi_{(f,g)}$ is differentiable and injective. The dynamics of $f : M \to M$ is mirrored using $\Phi_{(f,g)}$ by:

$$\Phi_{(f,g)}(f(x)) = (g(f(x)), g(f^2(x)), ..., g(f^{2m+1}(x))).$$

The theorem applies for function at least twice differentiable and continuous second derivative. If $f$ possess an attractor then this should show up in $\mathbb{R}^{2m+1}$ using the embedding.

## 2.12 Periodogram

Let $0 \leq \omega \leq \pi$ the periodogram ordinates are defined to be

$$I(\omega) = \left[\left(\sum_{t=1}^{n} y_t \cos(\omega t)\right)^2 + \left(\sum_{t=1}^{n} y_t \sin(\omega t)\right)^2\right]$$

The periodogram is plot of $I(\omega)$ against $\omega$.

There are some common conventions that are usually applied. The term $I(0)$ is usually omitted as it reduces to the sample mean and is of little direct interest in this scenario. Even though the term $I$ is defined for all $\omega$ it is usual to plot $I(\omega)$ for the positive Fourier frequencies only.

# 3 Modelling Methods

In this section the topic of modelling data will be discussed. Also some methods that have been previously been used to model systems will be discussed and the development will be shown. This will help when it comes to discussing the methods that have been chosen, as the previous methods will serve as a comparison.

## 3.1 Stochastic Models

Generally the most prominent stochastic models are the autoregressive model and the moving average model. These are of relevance as they find their roots in classical analysis. In general a stochastic model will consist of a process acting on a series of independent noise inputs and the past values of the signal itself. In this sense the methods can also be referred to as filters.

### 3.1.1 Linear Filters

In this set of stochastic models the processed used to filter the noise and signal inputs will be linear. When modelling such a system all the parameters are estimated from the outputs.

If the estimated power spectrum of the time series to be used has no prominent peaks it is in the form of coloured noise. For this type of data the moving average model is a good choice. It is a filter on a series of Gaussian white noise inputs $\eta_n$.

$$x_n = \sum_{j=0}^{M_{MA}} b_j \eta_{n-j}, \tag{1}$$

where $\langle \eta_n, \eta_m \rangle = \sigma^2 \delta_{nm}$[1] and $\langle \eta \rangle = 0$. Usually, $b_0$ is fixed and the number $M_{MA}$ of adjustable parameters is called the order of the process. Note that $x_n$ is also a Gaussian random variable with zero mean[2]. Sometimes, such a model is called a finite impulse response filter, since the signal dies after $M_{MA}$ steps, if the input is a single pulse, $\eta_0 = 1$ and $\eta_n = 0$ for $n \neq 0$.

If the power spectrum of a time series is dominated by peaks at distinct frequencies the autoregressive model is more appropriate. In a autoregressive model, the present outcome is a linear combination of the signal in the past (with finite memory), plus additive noise

$$x_n = \sum_{j=1}^{M_{AR}} a_j x_{n-j} + \eta_n. \tag{2}$$

This describes an autoregressive model of order $M_{AR}$, where $\eta_n$ is white Gaussian noise as above. Thus again, $x_n$ is a Gaussian random variable.

In principle, all Gaussian linear stochastic processes can be modelled with arbitrary accuracy by either if the two approaches. The number of parameters (i.e. the order of the model), however, might be extremely large, or it could even be infinite, if modelling a noisy harmonic oscillation by a moving average process.

The autoregressive moving average process, is as the name would suggest the combination of the autoregressive and moving average models. This yields a power spectrum with both pole and a polynomial background. The properties of all three models are only well understood if the increments, or inputs $\eta$, are Gaussian white noise. In general real world data is not Gaussian distributed. If these processes are used for real world data, usually it is assumed that a nonlinear transformation distorts the output of the Gaussian random process and changes the distortion to the one observed. Such nonlinearities are called static, in contrast to those of a dynamics of a system. These static nonlinearities conserve the property of time reversal invariance, which is also characteristic for linear stochastic processes. Before fitting an autoregressive moving average model to such data, the distribution should be rendered Gaussian by empirically inverting the nonlinear transformation.

---

[1] $\sigma$ being the standard deviation

[2] This can be achieved by subtracting the mean from a time series first.

Still the optimal order of these models remains unknown. Since the reproduction of data is better the more parameters the model contains, it is necessary to employ a standard which limits the number of coefficients in order to prevent overfitting.

The distinction of autoregressive, moving average and the combination models and the determination of the corresponding coefficients form a time series are only relevant if we are interested in the model *per se*. If the purpose of the modelling is to make forecasts of future values, recall that the noise inputs $\eta_n$ are not known to us and must be averaged over. This leaves only the autoregressive part as a predictive model.

### *3.1.2 Nonlinear filters*

Nonlinearity can be introduced to generalise an autoregressive models to produce such methods as the threshold autoregressive models [3]. A single threshold autoregressive model consists of a collection of the usual autoregressive models where each single one is valid only in a certain domain of the delay vector space (separated by the 'thresholds'). The construction of the model is performed by dividing the reconstructed state space into patches, and determining the coefficients of each single autoregressive model as usual, using only data points in the corresponding patch. Thus the threshold autoregressive models are piecewise linear models.

### 3.1.3 Markov models

Autoregressive models are a special case of Markov models. Markov models rely on the notion of a state space, although 'state' cannot be interpreted in the sense of determinism. In a Markov model of order *m* the probability of finding a signal at time *n* in some interval *i* depends only on the last *m* time steps. These last *m* time steps define the state of the system. On a scalar time series, a Markov model relies by definition on time delay embedding. Then a discrete Markov model is defined by specifying all transition probabilities from one cell of the state space to another. For delay coordinates one usually induces a partition in the state space by a partition on the interval, and then defines transition probabilities from the cells to the intervals, $p(i|j_1,...,j_m)$. For each cell $(j_1,...,j_m)$ one generally has several nonzero transition probabilities: otherwise the model would be deterministic. An autoregressive model is a Markov model of order $M_{AR}$ with continuous state space.

## 3.2 LOCAL METHODS

Local method consist of local neighbourhood to neighbourhood maps in the reconstructed state space. These are conceptually simpler than global methods, but as with stochastic models, depending on the purpose they can require a large numerical effort. In general local methods are useful for short predictions when new samples are being added to the time series. There are two main classes of local methods, those applying neighbour samples directly in the prediction, and those fitting a function locally to the neighbours basing the prediction on the estimated function.

### 3.2.1 Neighbour samples methods

The simplest way to predict a future point, $x_{k+1}$, of the time series from neighbour samples is to identify the nearest neighbour to the original point, $x_k$, in the embedding space $\mathbb{R}^m$. The nearest neighbour to $x_k$ is denoted by $x_{k(1)}$, and the next sample $x_{k(1)+1} = f(x_{k(1)})$ is then known from the time series, and can be used as the predictor. This is also termed the "analog method". An improvement is to take the *s* nearest neighbours and use the average of their state mappings as the predictor.

### 3.2.2 Fitting methods

An example of a fitting method is the local linear model. The deterministic part of this model is an autoregressive model of order *m*, such that the prediction is a linear superposition of the *m* last observables. The noise term responsible for the stochastic nature of the autoregressive mode is now interpreted slightly differently. It is now assumes that deterministic dynamics are being dealt with, but the observed data is contaminated by measurement noise.

## 3.3 GLOBAL METHODS

The idea of global modelling is straightforward. An appropriate functional form for *F* has to be chosen so that it is flexible enough to model the true function on the whole attractor. A very popular strategy is to take *F* to be a linear superposition of basis functions, $F = \sum_{i=1}^{k} a_i \Phi_i$. The *k* basis functions $\Phi_i$ are kept fixed during the fit and only the coefficients $a_i$ are varied.

### 3.3.1 Polynomials

Global linear models can be regarded as the first approximation in a Taylor series expansion of $F(s)$. The generalisation of this is to use a polynomial. Since *F* can act on *m*-dimensional space, it has to be a multivariate polynomial, which for order *l* has $k = (m+l)!/m!l!$ independent coefficients. These coefficients can usually be readily determined despite their large number as long as the data does not contain to much noise. Polynomials have the advantage that they are very familiar, which allows the result to be understood. The *k* coefficients can be determined by the inversion of a single ($k \times k$) matrix. Sometimes, however, prediction functions obtained by polynomial fits give rise to unstable dynamics under iteration since the polynomial basis functions diverge for large arguments and trajectories may escape to infinity.

### 3.3.2 Radical basis functions

An alternative to this basis of polynomials is to use a basis of radicals. These model are called, unsurprisingly, radical basis functions. One defines a scalar function $\Phi(r)$ with only positive arguments $r$. Additionally, one has to select $k$ centres $y_i$ on the attractor. Including the constant function as the zeroth basis function, $F$ reads

$$F(x) = x_0 + \sum_{i=1}^{k} a_i \Phi(|x - y_i|)$$

Typical basis functions $\Phi$ are bell shaped with a maximum at $r=0$ and a rapid decay towards zero with increasing $r$. Also increasing functions or even singular once can be used. The function $F$ is modelled by adjusting the coefficients $x_i$ of the functions $\Phi$. If the centres $y_i$ are reasonably well distributed on the attractor the superposition gives a well behaved surface.

Some consideration has to be taken when choosing the number and positions of the centres $y_i$. The centres and other parameters (width, etc.) of the basis functions $\Phi(r)$ are kept fixed. This makes the determination of the coefficients $a_i$ a linear problem. The typical width of the function $\Phi$ can be optimised by systematically testing several values, since for every fixed value the numerical fit procedure is very fast. Fitting then becomes nonlinear and thus is much harder when the centres are also to be optimised.

### 3.3.3 Neural Networks

Neural networks are analytic techniques modelled after the (hypothesised) processes of learning in the cognitive system and the neurological functions of the brain and capable of predicting new observations (on specific variables) from other observations (on the same or other variables) after executing a process of so-called learning from existing data. In fact, the networks used in this section are again nothing but nonlinear global models of the dynamics, the structure of which can be understood quite easily.

One particular class which has been used for time series modelling comprise of the so called feed forward networks with one hidden layer. This means that there is one layer of input units, one of neurons and one of output units. In the case of scalar predictions the latter is a single unit, just summing up the output of the neurons. The input layer consists of *m* units if the working space is *m*-dimensional embedding space, and it does nothing but distribute the *m* components of the delay vectors to the neurons. Since the whole structure of the network is inspired by the nervous system, the function of a neuron is usually a smoothed step function, a sigmoid function such as $\Phi = 1/(1+\exp(\mathbf{b}\mathbf{x}-c))$. The whole network is thus nothing other than the function

$$F(x) = \sum_{i=1}^{k} \frac{a_i}{1+\exp(\mathbf{b}_i \mathbf{x} + c_i)},$$

where the parameters $a_i$, $\mathbf{b}_i$ and $c_i$ have to be determined by a fit. The minimisation problem is nonlinear in the parameters, so in order to avoid overfitting the network should not be chosen to be too large. The standard way to determine the parameter is, in this context, called training. The most popular training algorithm is the error backpropagation. This is a gradient descent method where a cost function is minimised by presenting all learn pairs $(x_{n+1}, \mathbf{x}_n)$ of the one step prediction error individually to the function *F*. Then parameters are iteratively modified.

## 4 MEASURE BASED RECONSTRUCTION

In the spirit of a case study the Measure Based Reconstruction method by Giona, Lentinti and Cimagalli, described in the paper [4] will be investigated in this section. This is once again a global method of modelling data. The method is interesting as it allows the problem of multiparameter optimisation to be eliminated. The basis of the method is the idea of reconstructing *f(x)* from the knowledge of the time series $\{x_i\}$ by means of the generalised Fourier expansion with respect to the system of polynomials orthonormal to the invariant measure $\mu_I$ associated to $\{x_i\}$ and ultimately to *f(x)*. This means that the coefficients *c* are evaluated directly from the time series in terms of the hierarchy of moments **M** and the hierarchy of functional moments $\Gamma$ associated to $\mu_I$. The consequence of the proposed method is that the expansion coefficients *c* are $L^2$[3] optimal with respect to the natural metrics $\rho_I$ and are independent. The measure based approach is valid for time series arising from ergodic dynamical systems for which $\mu_I$ exists. This assumption is not so restrictive as most dynamical systems of physical interest can be considered as ergodic ones.

In this section the one dimensional, multidimensional discrete time models and the continuous time models will be presented. These methods, apart from the continuous time model, will be investigated using numerical experiments in the next section. This will test the methods for their ease of use and goodness of fittings.

## 4.1 ONE DIMENSIONAL DISCRETE TIME METHOD.

Consider a discrete dynamical system

$$x_{n+1} = f(x_n), \quad x_n \in \mathbb{R} \tag{3}$$

generating the time series $\{x_i\}$. The basic idea of the measure based approach lies in approximating *f* in terms of the polynomial system $\Pi = \{\pi^{(i)}(x)\}$, where $\pi^{(i)}(x)$ is a polynomial in *x* of degree *i*, associated to the invariant measure $\mu_I$ of (3)

$$\langle \pi^{(i)}, \pi^{(j)} \rangle_I = \int_\varrho \pi^{(i)}(x) \pi^{(j)}(x) d\mu_I(x) = \delta_{ij} \tag{4}$$

where $\varrho$ is the limit set on which the ergodic measure $\mu_I$ is concentrated and $\delta_{ij}$ is the Kronecker tensor. The system $\Pi$ can be called the *natural polynomial system* associated to $\mu_I$ because it is uniquely defined by the invariant measure $\mu_I$.

The $n^{\text{th}}$ order polynomial approximation $f_n(x)$ of *f* reads

$$f_n(x) = \sum_{k=0}^{n} c^{(k)} \pi^{(k)}(x) \tag{5}$$

---

[3] Continuous twice differentiable.

in which from ergodic theorem, the expansion coefficients $c^{(k)}$ are given by

$$c^{(k)} = \int_\varrho f(x)\pi^{(k)}(x)d\mu_I(x)$$
$$= \lim_{N\to\infty} \frac{1}{N} \sum_{i=1}^{N} x_{i+1}\pi^{(k)}(x_i). \quad (6)$$

The polynomial system $\Pi$ can be obtained from the knowledge of the hierarchy of moments $M = \{M_k\}$:

$$M_k = \int_\varrho x^k d\mu_I(x) = \lim_{N\to\infty} \frac{1}{N} \sum_{i=1}^{N} x_i^k \quad (7)$$

by making use of the Gram-Schmidt orthogonalisation method. The practical calculation of the $\pi^{(i)}$'s can be performed as an iterative procedure starting from

$$\pi^{(0)}(x) = 1$$

which is a consequence of (4). If $\Pi_{k-1} = (\pi^{(0)}(x), \ldots, \pi^{(k-1)}(x))$ is known and

$$\pi^{(h)}(x) = \sum_{j=0}^{h} a_j^{(h)} x_j \quad (8)$$

a polynomial of order $k$, orthonormal to $\Pi_{k-1}$, is given by

$$P_k(x) = x^k + \sum_{j=0}^{k-1} \beta_j^{(k)} \pi^{(j)}(x) \quad (9)$$

where

$$\beta_j^{(k)} = -\sum_{i=0}^{j} a_i^{(j)} M_{k+i}, \; j = 0, \ldots, k-1 \;. \quad (10)$$

The coefficients $a_i^{(k)}$ of $\pi^{(k)}(x)$ are given by

$$a_i^{(k)} = \frac{\left[\sum_{k=0}^{k-1} \beta_h^{(k)} a_i^{(h)}\right]}{\sqrt{N_k}}, \; i \neq k$$
$$= \frac{1}{\sqrt{N_k}}, \; i = k \quad (11)$$

Where the normalisation constant $N_k$ reads

$$N_k = \langle P_k, P_k \rangle_I = M_{2k} - \sum_{h=0}^{k-1} (\beta_h^{(k)})^2 \geq 0. \quad (12)$$

It can be seen from (6) that the coefficient $c^{(k)}$ can be expresses in terms of the hierarchy of *functional moments* $\Gamma = \{\Gamma^{(j)}\}$, defined as

$$\Gamma^{(j)} = \int_\varrho f(x) x^j d\mu_I(x)$$
$$= \lim_{N \to \infty} \frac{1}{N} \sum_{i=1}^{N} x_{i+1} x_i^j, \quad j \in \mathbb{N}^+, \tag{13}$$

obtaining

$$c^{(k)} = \sum_{j=0}^{k} a_j^{(k)} \Gamma^{(j)}. \tag{14}$$

As can be seen from the above analysis the measure based reconstruction is based exclusively on the evaluation of the ergodic sums (7) and (13). Therefore it does not need any kind of parameter fitting. This method has been programmed in to the MBRecon program in C, which is listed at the end of the project.

## 4.2 Extension to Multidimensional Dynamics

As one dimensional method introduces the basic principles of measure based reconstruction, this section will now look at the extension to multidimensional dynamics. This should cause no great theoretical problems, however it is convenient to develop the analysis at least for two-dimensional dynamics. This is so the structure of the polynomial system $\Pi$ can be clarified.

In two dimensions, $(x_1, x_2)$, $\Pi = \{\pi^{(i,j)}(x_1, x_2)\}$ is given by

$$\pi^{(i,j)}(x_1, x_2) = \sum_{k=0}^{i-1} \sum_{h=0}^{k} a_{hk}^{(i,j)} x_1^{k-h} x_2^h + \sum_{k=0}^{j} \eta_k^{(i,j)} x_1^{i-k} x_2^k \tag{15}$$

and the maps $f_s(x_1, x_2)$, $s = 1, 2$ are approximated by

$$f_s(x_1, x_2) = \sum_{i=0}^{n} \sum_{j=0}^{i} c_s^{(i,j)} \pi^{(i,j)}(x_1, x_2), \quad s = 1, 2. \tag{16}$$

Of course $\eta_0^{(0,0)} = 1$, as from the one dimensional method we would expect $\pi^{(0,0)}$, the zeroth ordered polynomial to be equal to 1.

The generic polynomial $P^{(i,j)}(x_1, x_2)$, orthogonal to the preceding polynomials $\pi^{(m,r)}$ of the system $\Pi$, $0 \leq m \leq i$, $0 \leq r < j$, can be written as

$$P^{(i,j)}(x_1, x_2) = \sum_{h=0}^{i-1} \sum_{k=0}^{h} A_{hk}^{(i,j)} \pi^{(h,k)}(x_1, x_2) + \sum_{k=0}^{j-1} B_k^{(i,j)} \pi^{(i,k)}(x_1, x_2) + x_1^{i-j} x_2^j$$

where the coefficients $A_{hk}^{(i,j)}$, $B_k^{(i,j)}$ read

$$\begin{aligned}
A_{hk}^{(i,j)} &= -\sum_{l=0}^{h-1} \sum_{m=0}^{l} a_{ml}^{(h,k)} M_{l-m+i-j,m+j} - \sum_{m=0}^{k} \eta_m^{(h,k)} M_{h-m+i-j,m+j}, \\
B_k^{(i,j)} &= -\sum_{l=0}^{i-1} \sum_{m=0}^{l} a_{ml}^{(i,k)} M_{l-m+i-j,m+j} - \sum_{m=0}^{k} \eta_m^{(i,k)} M_{2i-m-j,m+j}
\end{aligned} \tag{17}$$

and the hierarchy of moments $\mathbf{M} = \{M_{i,j}\}$ is given by

$$M_{i,j} = \int_\varrho x_1^i x_2^j d\mu_I(x_1, x_2) = \lim_{N \to \infty} \frac{1}{N} \sum_{k=1}^{N} x_{1,h}^i x_{2,h}^j. \tag{18}$$

After some algebraic manipulations, we obtain for the $a^{(i,j)}$'s and $\eta^{(i,j)}$'s the expressions

$$\begin{aligned}
a_{hk}^{(i,j)} &= \frac{\sum_{r=k+1}^{i-1} \sum_{q=0}^{r} A_{rq}^{(i,j)} a_{hk}^{(r,q)} + \sum_{q=h}^{k} A_{kq}^{(i,j)} \eta_k^{(k,q)} + \sum_{q=0}^{j-1} B_q^{(i,j)} a_{hk}^{(i,q)}}{\sqrt{N_{i,j}}}, \\
\eta_k^{(i,j)} &= \frac{\sum_{q=k}^{j-1} B_q^{(i,j)} \eta_k^{(i,q)}}{\sqrt{N_{i,j}}}, \quad k \neq j \\
\eta_k^{(i,j)} &= \frac{1}{\sqrt{N_{i,j}}}, \quad k = j
\end{aligned} \tag{19}$$

where

$$N_{i,j} = \langle P^{(i,j)}, P^{(i,j)} \rangle_I = M_{2(i-j),2j} - \sum_{h=0}^{i-1} \sum_{k=0}^{h} (A_{hk}^{(i,j)})^2 - \sum_{k=0}^{j-1} (B_k^{(i,j)})^2. \tag{20}$$

As in the one-dimensional case the coefficients $c_s^{(i,j)}$ can be obtained from the hierarchy of moments $\Gamma = \{\Gamma_s^{(h,k)}\}$:

$$\Gamma_s^{(h,k)} = \lim_{N \to \infty} \frac{1}{N} \sum_{i=1}^{N} x_{s,i+1} x_{1,i}^h x_{2,i}^k, \quad s = 1, 2, \ h, k \in \mathbb{N}^+ \tag{21}$$

which gives

$$c_s^{(i,j)} = \sum_{k=0}^{i-1} \sum_{h=0}^{k} a_{hk}^{(i,j)} \Gamma_s^{(k-h,k)} + \sum_{k=0}^{j} \eta_k^{(i,j)} \Gamma_s^{(i-k,k)}. \tag{22}$$

Despite the apparent formal complexity of the measure based reconstruction, this method presents great advantages, as discussed in the one dimensional case, with respect to the other methods of functional reconstruction. The generalised Fourier coefficients *c* can be directly evaluated starting from the hierarchies of the moments *M* and Γ. Furthermore the Fourier expansion ensures that the coefficients $c_s^{(i,j)}$ are totally independent of each other.

## 4.3 CONTINUOUS TIME METHOD

It is easy to see that the measure based approach can be applied with some modification even in the continuous case. The fundamental structure of measure based reconstruction remains, however, unchanged in the transition from discrete to continuous dynamical systems. In the continuous case a vector field *f*:

$$\dot{x} = f(x)$$

has to be estimated from a time series $\{x_i\}$.

Considering three-dimensional phase space $x = (x_1, x_2, x_3)$, and using the previously discussed technique, the natural polynomial system $\Pi$ can be expressed in terms of the hierarchy of moments $\{M_{i,j,k}\}$:

$$M_{i,j,k} = \int_\varrho x_1^i x_2^j x_3^k d\mu_l(x). \qquad (23)$$

However, the expansion coefficients $c_s^{(i,j,k)}$:

$$f_s(x_1, x_2, x_3) = \sum_{i,j,k} c_s^{(i,j,k)} \pi^{(i,j,k)}(x_1, x_2, x_3), \ s = 1, 2, 3$$

have to be estimated from the infinitesimal functional moments $\Theta_s^{(i,j,k)}$:

$$\Theta_s^{(i,j,k)} = \lim_{N \to \infty} \sum_{m=1}^{N} \left[ \frac{x_{s,m+1} - x_{s,m}}{\Delta t} \right] x_{1,m}^i x_{2,m}^j x_{3,m}^k, \ s = 1, 2, 3 \qquad (24)$$

which, in practise, depends upon the sampling time $\Delta t$.

The mathematical formulation in this case is completely equivalent to that previously discussed and is therefore omitted.

# 5 NUMERICAL EXPERIMENTS

At this point in the project enough theory has been covered to perform some numerical experiments using the Measure Based Reconstruction Method. The variables will be sequentially fixed so the affect of the prediction length and other measure of the data by the method can be evaluated.

For these experiments the programs listed at the end of the project will be used. These programs are written so that the interested reader can perform other investigations using the measure based method. For the execution of the programs a Pentium III 1 GHz and a Pentium 233 MHz have been used and no noticeable difference has been noted in the execution time. However, the MATLAB package has been used to produce the graphical elements of the section, and this I advise should be run on the Pentium III.

## 5.1 METHOD

In our experiments, apart from those performed on real data, the data will be generated as is needed by the Quadratic program. The noise will be added to these synthetic data sets using the addnoise package from the TISEAN package [24]. The data generation program is only used as it could be adapted to fit the generation method that was needed.

The maps that will be used to generate the data are the following:
1. Exponential quadratic, $x_{n+1} = \mu x_n (1-x_n) \exp(-kx)$ with $\mu = 10$ and $k = 2.51705$. Initial point 0.123456789
2. Quadratic, $x_{n+1} = \mu x_n (1-x_n)$ with $\mu = 3.8$. Initial point 0.123456789.
3. Hénon $x$-variable, $x_{n+1} = 1 + bx_{n-1} - ax_n^2$ with $a = 1.4$ and $b = 0.3$ Initial points 0.12 and 0.22.

The first two of these maps will produce trajectory data and so the method should be able to predict to a good length (more than 5) within a certain error value, $\varepsilon$. The Hénon data set will not be trajectory data as it is a two dimensional attractor, but it will serve to prove that the measure based reconstruction can work on non trajectory data.

## 5.2 Proving the Paper

To start with the one dimensional method will be examined. The first experiment was try and replicate the results that have been given in the paper [4]. The fitting of these initial functions will be examined to see if the residuals contain any pattern. If they do then it means that some of the information that should be preserved by the model is not being carried though.

Below is the graph that is trying to be replicated. As the prediction will be measured from the end of the data set the prediction length will be taken as $T(\varepsilon)$. Also note that the graphic below does not have a uniform scale and so the graphs which result from the experiment will possibly look slightly different.

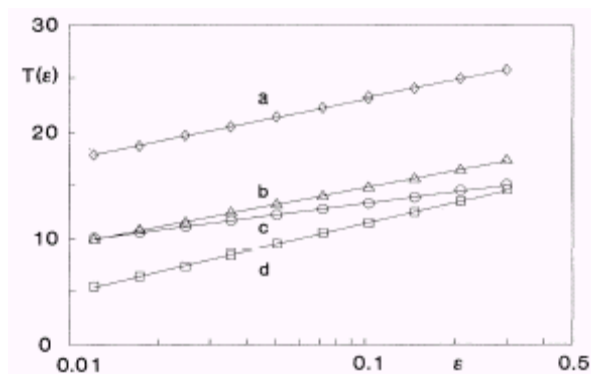

FIG. 1. $T(\epsilon)$ vs $\epsilon$ for one-dimensional and two-dimensional maps: (a) Hénon map, $n=2$, $N=50\,000$; (b) exponential quadratic map with $p=10$, $k=2.517\,05$, $n=10$, $N=30\,000$; (c) quadratic map with $p=3.8$, $n=2$, $N=30\,000$; (d) Ikeda map $n=5$; $N=50\,000$.

### 5.2.1 Prediction

For this experiment the prediction is that the results which are obtained will follow those from the paper, but may be smaller. This is because the claims in the papers, as we can see from the graphs above seem a little ambitious. The normal prediction length that would be expected from a method would be about 5 points.

*5.2.2  The Results*

The graph of the resultant prediction lengths of the reconstruction is shown below:

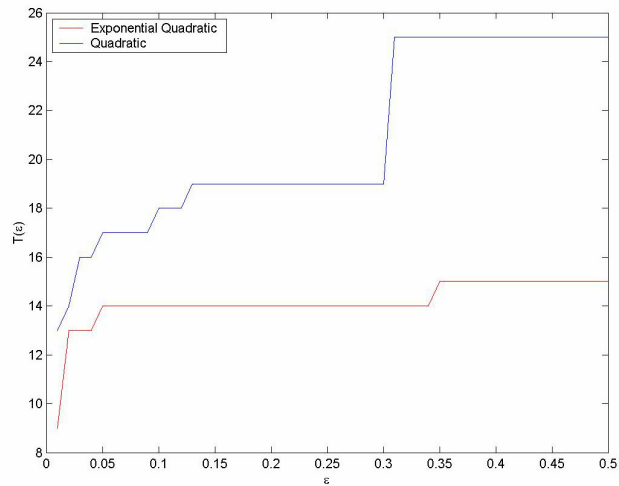

As we can see the initially that the results seem to be inverse of the results from the paper. The graphs seem to be generally increasing in length as the tolerance is greater, although there are also plateaux in both graphs suggesting that the prediction remains at the same accuracy. The magnitudes of the prediction lengths are also higher than those in the paper which could be expected as the data set used is most likely different. The reason for the inverse is that the quadratic data is substantially better predicted than it was in the paper. This may be due to the method used to compute the prediction length of the data set. In the paper an average over the whole data set is taken to determine the mean prediction length. This means that the prediction length may fluctuate and so produce a lower value. In our experiment the prediction length at the end of the training set.

Now as the quadratic data seems to have different qualities to the quadratic result given in the paper, the data sets will be analysed individually. We shall start with the exponential quadratic data as it seem to follow the pattern from the paper. To start with let us look at the prediction length graph again.

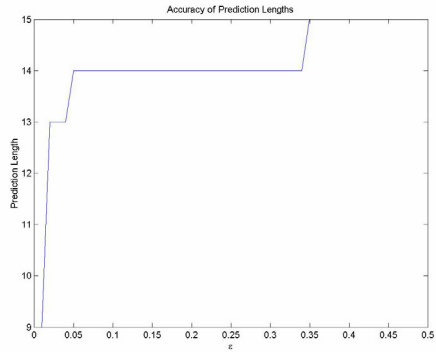 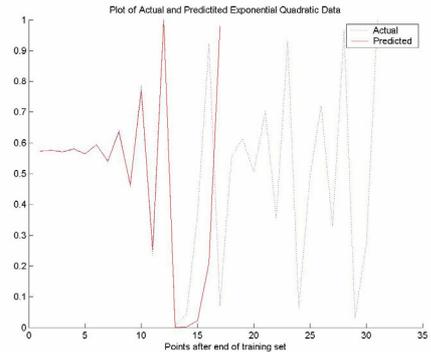

As we can see the prediction length achieves a maximum of 15 points at the tolerance level of 0.35. This would relate to a 35% noise level and this is a little high for most prediction applications. In the more sensible tolerance levels between 0 and 0.05, the prediction length is at most 13 points and at least 9. This length is unusual in most prediction methods, so we can say that it is good for prediction.

To extend our study of this data set and the fitting we shall look at the periodograms of both time series. This will be done on the same axis so that we can easily seem and differences.

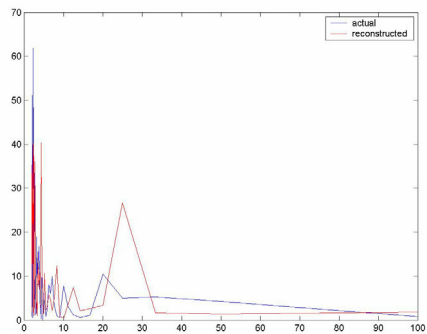 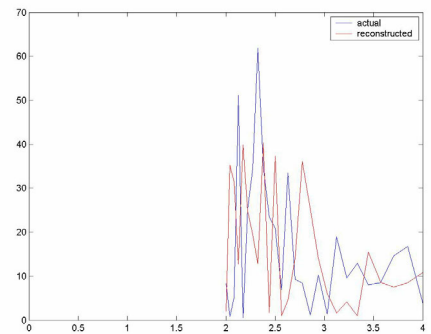

As we from the graph on the left most of the interesting peaks are in the area between 2 and 4. The reconstructed data however have a peak in the 20-30 interval that seems quite large for it's position. This shows that the reconstruction's accuracy was not consistent all of the time with the actual data. The graph on the right show the interesting region zoomed in. As we can see[4] the peak of the actual data is at 2.33 with a power of 61.88 while the reconstructed data has its peak at 2.38 with power 40.37. Also the reconstructed data appear to have a second peak of about the same height suggesting another periodic component. This difference in the position of the peak and the reconstructed data can be expect as the reconstruction only predicts 13 points accurately. Overall the power spectrum of the reconstructed data is lower and contains more peaks, suggesting more periodic components but generally it follows that of the actual data.

---

[4] Well you might be able to see but it is easier reading the MATLAB output

The lag plot of the data should also be looked at to see if the residuals have any structure. This plot is shown below

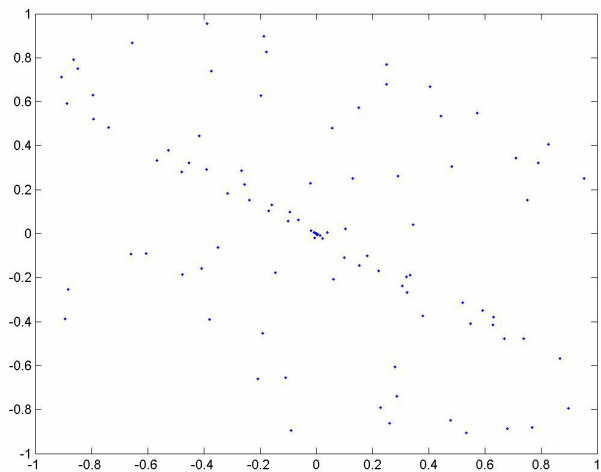

There is a small amount of structure to the residuals from the fitting, suggesting a correlation between the points. However if we were to take out the first 15 points the structure would probably be lost. Generally the noise appears to have no pattern and so could be considered as the noise for the repeated computations in the program and rounding errors.

Now we shall move on to look at the quadratic data, and try and understand why it might be a better prediction than had been achieved in the paper. Let us, as before look at the prediction graphs.

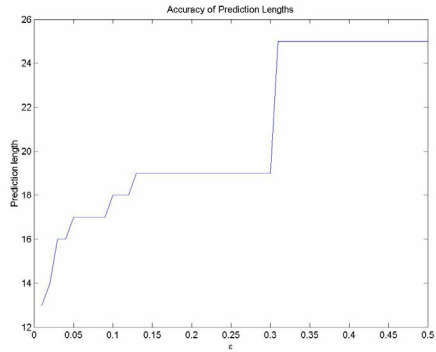 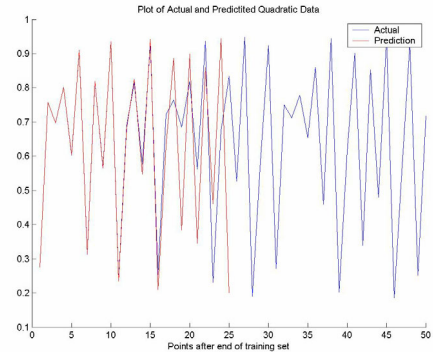

Once again we see the plateaux in the prediction length graph and that it achieves a prediction length of 25 points at a tolerance of just over 0.3. This would seem better than the paper's prediction level, but once again it must be pointed out that the paper is using a mean prediction level. Back at our sensible levels for the tolerance the prediction length appears to be between 13 and 16 points, once again beating most prediction methods commonly used. We can see from the sequence graph on the right that the reconstructed data follows the actual data quite closely and then the amplitudes of the peaks become higher than those witness in the real data, although they are at the same points.

Below is the periodogram for this data:

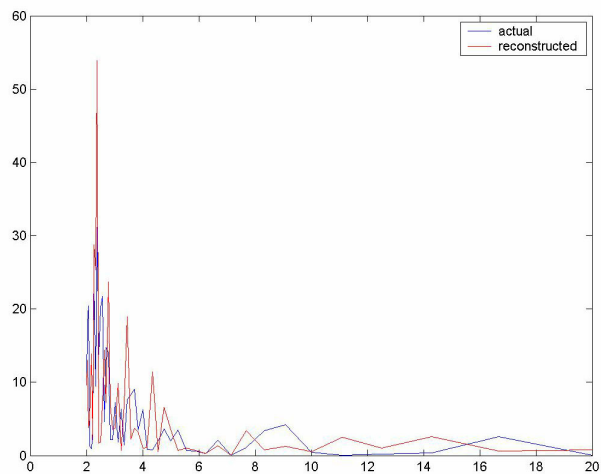

As we can see from the periodogram, to a certain level of accuracy the data sets had the period value for their spike. However this time it is the reconstructed data that has the greater power. The value for the data sets are that they have their peaks at period 2.3810 and the actual data has a power of 33.18 while the reconstructed data has a power of 53.91. Once again it can be noted that the reconstructed data set's spectrum is contains more peaks but most of them see to be where the actual data set's one are.

The lag plot of the residuals of the data is shown below:

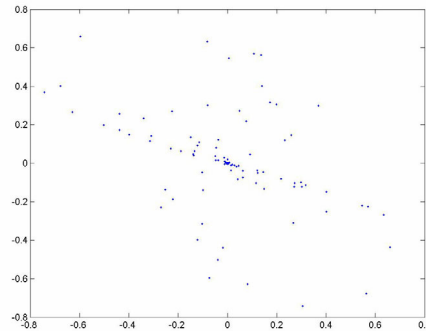

Once again there appears to be a pattern of a line towards the middle of the plot but this is probably the predictive points who are close to the actual data set. If these are removed again it could quite easily be seen that the plot looks similar to one of noise. This suggest that the structure of the real data set is remaining in the model.

*5.2.3   Conclusion*

The results of the experiment seem to general follow those given in the paper. The magnitude of the prediction lengths is slightly higher but that could just be the difference in methodology. We have also seen that the modelling technique works in the one dimensional case quite well, and seems to be easy to apply with short computation times. The modelling seems to be close.

## 5.3   NOW FOR TWO DIMENSIONS.

For our second experiment the aim will to be to try a similar process to the previous experiment. Once again the claims from the paper will be tested against those which are obtained from performing the reconstruction ourselves. Again the residuals will be investigated for any pattern, to see if the model is carrying through the information of the system.

### 5.3.1 Prediction

As the last experiment showed the claims that have been given in the paper are quite close to the experimental results that can be obtained independently. Knowing this, the prediction is that the experiment will produce values which are close in style to those in the paper although they might has slight differences in magnitude.

### 5.3.2 A slight problem

While trying to carry out the experiment with the program MBRecon2, listed at the end of the project, a small problem was noted that resulted in the results of the experiment being void. In this section the reasons for the problem will be investigated and discussed. It does however mean that for the rest of the experiments, they can only be performed in one dimension.

After debugging the program for the two dimensional case, the problem seemed to be caused by the $N_{2,0}$ being negative. As this number has to be square rooted the fact it is negative means that it produces a complex result. Not something that was expected when dealing with real data. To determine whether the program was at fault or the algorithm the $N_{2,0}$ value was worked out from the formulas and then calculated when the final equation was found.

To work out the value of $N_{2,0}$ the following formula is used

$$N_{2,0} = M_{4,0} - \sum_{h=0}^{1}\sum_{k=0}^{h}\left(A_{hk}^{(2,0)}\right)^2 \qquad (25)$$
$$= M_{4,0} - \left(A_{00}^{(2,0)}\right)^2 - \left(A_{10}^{(2,0)}\right)^2 - \left(A_{11}^{(2,0)}\right)^2$$

Now the following quantities have been derived from the equation is Section 4:

$$B_0^{(1,0)} = \frac{(M_{1,0})^2 \eta_0^{(0,0)}}{\sqrt{N_{1,0}}} - \frac{M_{2,0}}{\sqrt{N_{1,0}}}$$

$$B_0^{(1,1)} = -a_{00}^{(1,0)} M_{0,1} - \eta_0^{(1,0)} M_{1,1}$$

where

$$a_{00}^{(1,0)} = \frac{A^{(1,0)} \eta_0^{(0,0)}}{\sqrt{N_{1,0}}}$$

$$\eta_0^{(0,0)} = 1$$

$$\eta_0^{(1,0)} = \frac{1}{\sqrt{N_{1,0}}}$$

and

$$A_{00}^{(1,0)} = -M_{1,0}$$
$$A_{00}^{(1,1)} = -M_{0,1}$$

So now to calculate the variables in (25) by substituting in the other values. First calculate $A_{00}^{(2,0)}$:

$$A_{00}^{(2,0)} = -\eta_0^{(0,0)} M_{2,0} = -M_{2,0}. \tag{26}$$

Now on to $A_{10}^{(2,0)}$:

$$\begin{aligned} A_{10}^{(2,0)} &= -a_{00}^{(1,0)} M_{2,0} - \eta_0^{(1,0)} M_{3,0} \\ &= -\frac{M_{1,0} \eta_0^{(0,0)}}{\sqrt{N_{1,0}}} - \frac{M_{3,0}}{\sqrt{N_{1,0}}} \\ &= -\frac{M_{1,0}}{\sqrt{N_{1,0}}} - \frac{M_{3,0}}{\sqrt{N_{1,0}}} \end{aligned} \tag{27}$$

The calculation for those variables from (1) were fairly straight forward. It is a different matter when it comes to $A_{11}^{(2,0)}$

$$A_{11}^{(2,0)} = -a_{00}^{(1,1)} M_{0,1} - \eta_0^{(1,1)} M_{1,1} - \eta_1^{(1,1)} M_{0,2} \tag{28}$$

As this expansion can get a bit tricky if done all at one, it shall be computed term by term. So

$$\begin{aligned} a_{00}^{(1,1)} M_{0,1} &= \frac{A_{00}^{(1,1)} M_{0,1} \eta_0^{(0,0)} - B_0^{(1,1)} a_{00}^{(1,0)} M_{0,1}}{\sqrt{N_{1,0}}} \\ &= \frac{-M_{0,1} - B_0^{(1,1)} a_{00}^{(1,0)} M_{0,1}}{\sqrt{N_{1,0}}} \end{aligned}$$

Now substituting for *B*,

$$a_{00}^{(1,1)} M_{0,1} = -\frac{M_{0,1}}{\sqrt{N_{1,0}}} + \frac{\left(M_{0,1} a_{00}^{(1,0)}\right)^2}{\sqrt{N_{1,0}}} + \frac{\eta_0^{(1,0)} a_{00}^{(1,0)} M_{0,1} M_{1,1}}{\sqrt{N_{1,0}}}$$

$$= -\frac{M_{0,1}}{\sqrt{N_{1,0}}} + \frac{(M_{0,1} M_{1,0})^2}{N_{1,0}^{\frac{3}{2}}} - \frac{M_{0,1} M_{1,1} M_{1,0}}{N_{1,0}^{\frac{3}{2}}}.$$

As this is in terms of the moments of the data it can not be simplified any further. Now onto the second term of the equation

$$\eta_0^{(1,1)} M_{1,1} = \frac{B_0^{(1,1)} \eta_0^{(1,0)} M_{1,1}}{\sqrt{N_{1,1}}}$$

$$= \frac{B_0^{(1,1)} M_{1,1}}{\sqrt{N_{1,0}} \sqrt{N_{1,1}}}$$

$$= \frac{M_{1,0} M_{0,1} M_{1,1}}{N_{1,0} \sqrt{N_{1,1}}} - \frac{(M_{1,1})^2}{N_{1,0} \sqrt{N_{1,1}}}.$$

Finally the last term is quite simple to compute

$$\eta_1^{(1,1)} M_{0,2} = \frac{M_{0,2}}{\sqrt{N_{1,1}}}.$$

Combining these three term we obtain the formula for $A_{11}^{(2,0)}$:

$$A_{11}^{(2,0)} = \frac{M_{0,1}}{\sqrt{N_{1,0}}} - \frac{(M_{0,1} M_{1,0})^2}{N_{1,0}^{\frac{3}{2}}} + \frac{M_{0,1} M_{1,1} M_{1,0}}{N_{1,0}^{\frac{3}{2}}} - \frac{M_{1,0} M_{0,1} M_{1,1}}{N_{1,0} \sqrt{N_{1,1}}}$$
$$+ \frac{(M_{1,1})^2}{N_{1,0} \sqrt{N_{1,1}}} - \frac{M_{0,2}}{\sqrt{N_{1,1}}} \qquad (3)$$

Finally when these equations for the *A*'s are substituted in to the equation for $N_{2,0}$ the result is obtained.

Using these formulas the value of $N_{2,0}$ is calculated for a Hénon data set, to see if the problem that was discovered in MBRecon2 is a problem with the program or the method. This was done using the proving m-file[5] in MATLAB. The result of the calculation is that the value of $N_{2,0}$ remains negative, and so produces a complex number when it square rooted. This means that all polynomial approximations using the measure based reconstruction method above the linear terms will return an error, or infinity.

---

[5] Listed in the program listings section.

*5.3.3 Conclusion*

Although the experiment was never executed in full due to the problem found in the method, it did produce interesting results. In this project when the method is described it is complete and has all the indices, something which is not true of its description in the original paper. This could just be a small mistake that happen during the transfer to the journal it was published in. The index that has been in inserted on $\eta$ in the calculation of the $a$'s is probably correct in the project, as it seem the most sensible and followed the other equations. Additionally to the work shown in this experiment the other indices we tested in the missing position, and apart from using *i* the result was identical. When *i* was used the $\eta$ value became zero, and so large amounts of the equations were equal to zero. This help in the sense that then the *N*'s were all positive, but when tested in the program it did not produce a good approximation.

## 5.4 PROBABLY WORTHLESS.

After the failure of the second experiment, and so the two dimensional method not working I had an idea. The Hénon data set is created using two formulae, one for each dimension of the data set, as follows:

$$x_{1,n+1} = 1 + x_{2,n} - 1.4x_{1,n}^2$$
$$x_{2,n+1} = 0.3x_{1,n}$$

As the second dimension's equation seem to use the previous value in the one dimension's time series, and then this value used to create a future value of the one dimension's time series, I wondered if they could be modelled using the one dimension method. This one dimensional reconstruction could the be sent through a delay coordinate embedding method and returned to two dimensions. The equation used to generate the data set would be

$$x_{n+1} = 1 + 0.3x_{n-1} - 1.4x_n^2.$$

After discussing this idea with my supervisor, I understood that it would probably have no meaning what so every and should not be included.

However after thinking about this I have decided to include it a different data type to that of the exponential and regular quadratic data, and so might show if the method can be implemented generally, as the paper suggested.

### 5.4.1  Prediction

As it is non standard data, and so probably doesn't fit the conditions for the method, the prediction length will be small. This is based on the method using a next point prediction based on the current point, while the Hénon data set requires both the current point and one past point to be used. However when the predicted points are embedded into two dimensional space, they will be topologically conjugate with the training set.

### 5.4.2  Results

The graph of the prediction lengths achieved for a given tolerance level is shown below

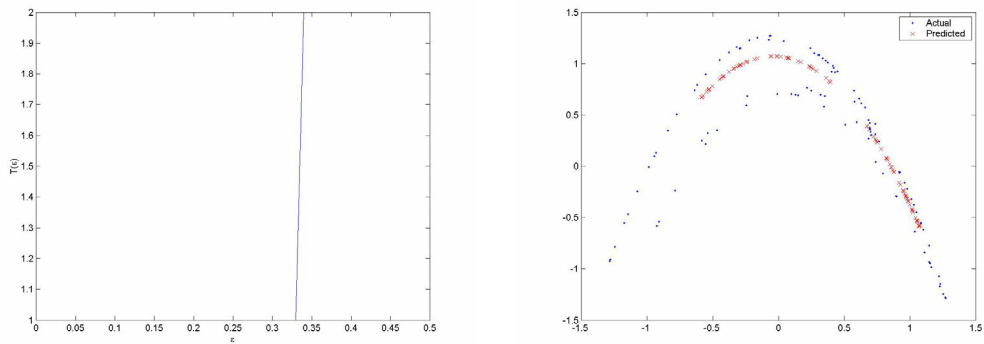

The prediction graph shows the maximum length to be two points with an tolerance level of about 0.34. This unlike the previous experiments means that the method is not very good. As the data is two dimensional the graph showing the predicted values against the actual values is two dimensional as well. The prediction points, when embedded, look slightly better when viewed in this context and seem to be at least topologically conjugate. The pattern of the Hénon attractor seems to have been kept in the points although the magnitude does seem to be a little reduced and so the curves of the attractor are closer together.

Now the periodogram is examined to determine whether the reconstructed data has similar properties to the training data. This is shown below:

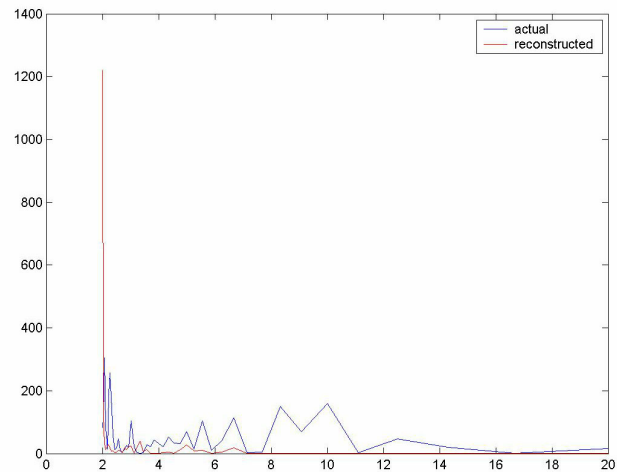

The predicted data and actual data seem to have similar peak locations, with the predicted data have its peak of power 1219 at period 2.00 while the actual data has a peak of 304.4 at period 2.08. Also the actual data seems to have a double spike of similar powers at the beginning and then seems to be more peaky than the predicted data. This suggests that the prediction is not very strong and does not capture all of the properties of the training set.

The lag plot of the residuals is

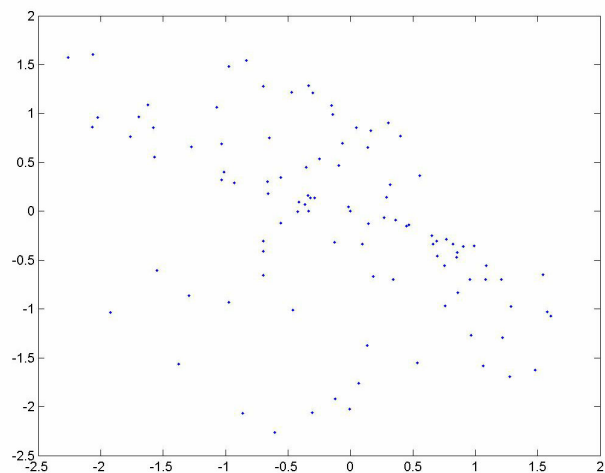

The lag plot appears to have more structure than in the previous experiments, once again suggesting that the fitting is not as good as we would like it to be. The structure however seems to be a wide shape and the residuals are not very densely packed. This suggests that there is still an element of noise in the prediction, possible again from the many computations and the rounding of the system.

*5.4.3    Conclusion*

This experiment has shown that the one dimensional reconstruction method can not be applied to data where the underlying process requires two points to predict a future value. It also showed that even though the predictions were not accurate most of the properties of the attractor were carried though to the predicted points. In general this showed that this approach to multidimensional data sets would not work, and that they should be reconstructed and then sent through the multidimensional version of the algorithm.

## 5.5    ROBUSTNESS TO NOISE

In this experiment the method will be tested for it robustness of prediction when noisy data is used. The method will be tested with different amounts of noise in the data to see the effects on the prediction length. If the method is robust to noise little effect should be noted. In the experiment the data used will be the quadratic data set, which will have Gaussian noise of the appropriate amplitude added to it using the `addnoise` function in the TISEAN package. The amplitude of the noise will be a percentage of the variance of the data, while the mean of the data will be zero. This means that the noise will be normally distributed and additive. The noise levels and amplitudes are shown in the table below.

| Noise level | Amplitude |
|---|---|
| 0.01 | 0.0406 |
| 0.05 | 0.2029 |
| 0.10 | 0.4058 |
| 0.15 | 0.6087 |
| 0.20 | 0.8116 |
| 0.30 | 1.2175 |
| 0.40 | 1.6233 |
| 0.50 | 2.0291 |

As this experiment is rather more qualitative than quantitative there will be less plots. Also there will be no proof that the predicted data will have all of the properties of the training sets as that is not being focussed on. Instead the difference in the increasing noise will be investigated by producing plots with multiple noise level on them.

### 5.5.1 Prediction

As the noise will some times be higher than our tolerance levels the prediction time will be reduced. This effect will become more prominent as the noise increases in size. Also the lengths of the plateaux will increase as the noise level increases, as the prediction points approaching the tolerance will be pushed over it by the noise. In general the prediction length will probably half.

### 5.5.2 Results

First the 1% to 15% noise level results will be discussed. The prediction length graph is shown below

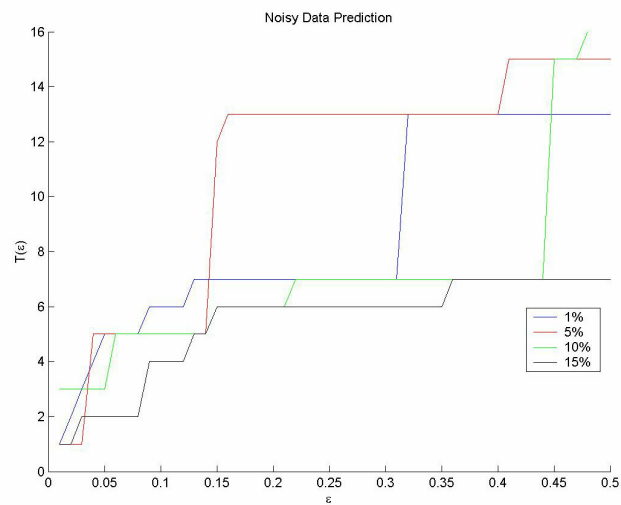

The prediction lengths of the noisy data look still quite good overall. In general they seem still to be quite high and so it would seem to be robust again noise up to 15%. Some interesting properties can also be noted. The first being that the data with 5%, and to a lesser extent 10%, seems to be able to be predicted better than the data with 1% noise. This is a strange result as it would be expected that the data with 1% noise would be the one with the best prediction lengths. It also appears that all the noise levels apart from the 10% can be predicted to the same length at a tolerance of 0.01. Once again the 10% noise level is predicted better than the rest of the data sets.

These properties that seem to be different to what would be expected could be explained by blind luck. It might be that the noise of the 10% noise data set was more beneficial than that on the 1%. In general though the prediction lengths seem to be sensible up to the level of tolerance, and then the noise is being accounted for. That is the prediction length is lower if there is a higher noise level.

Now the discussion moves onto a greater range of noise, to see if the pattern continues up to quite high levels of noise. The graph for this discussion is shown below.

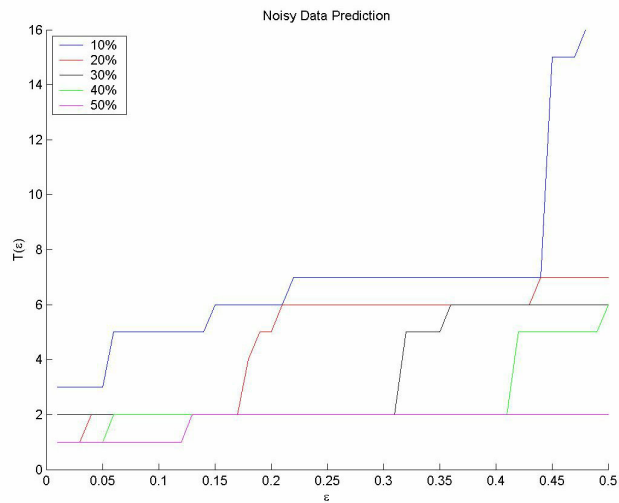

Again it can be noted that the prediction lengths get smaller as the noise gets larger. Also the decay of the prediction length appears to be quite fast, although the 10% value is still a little ambiguous in its values. The length of the plateaux seem also to get longer as the noise level increases. In fact by the time the noise reaches 50%[6] there seems to be only a short prediction that does not improve much as the tolerance is larger. This suggest that the method is quite robust to noise, as the length is still quite good with 30% noise.

---

[6] An artificially high value.

### 5.5.3 Conclusion

The experiment has shown that the method is very robust to noise. This is useful as it is quite unusual for real world data to be non noisy. The only possible exception to this is Stock market prices, but they still have external governing forces, and so may behave in a similar way. The one thing that should be noted is that the method does not get rid of the noise, it just tries to model the underlying process of the data. This means that it will also try to model the generation function for the noise. As noise is by definition random, it will be a bad approximation of this function as there is no connection between successive points. Also as the noise level gets higher the more the underlying process will be ignored for the process generating the noise.

## 5.6 CONCLUSION

From our experiments we have proved that the results of the paper can be replicated and that the method is robust to noise. Also we have seen that it cannot model data which is not in its correct dimension and that the two dimensional model might have a problem. The experiments were generally successful and a lot can be found out about a method like this. They have also shown that the analysis of real data would be possible if it was on the right form. It also seems that it can be quickly applied to most situation without much extra effort and is generally fast if run in a language similar to C.

Further experiments can be carried out on the method to find out more. These would include trying different sized data sets and seeing whether this effected the prediction length and changing the order of the approximation. These were not done in this project as the main aim of the experiments were to give a guide to those that should be carried out on modelling methods to evaluate their use. Although the size of the data set and the changing of the order of the approximation are important experiments, from similar methods the results can be predicted, and the answers would seem common sense as well. In general if the training set is bigger the model is better as it will be able to learn more about the behaviour of the system. If the order is increased the model will be closer to that of the data, but it will be at the cost of the computations required. When considering these experiments the actual advantage of making the quantity larger must be considered, as the extra computations may outweigh the advantage.

# 6 ALTERNATIVE METHODS

The Measure Based Reconstruction method, as has been seen, is a good modelling method as the coefficients of the approximation can be worked out from sums of the data. The method relies on the superposition of orthogonal polynomials to approximated to the function being modelled. This is unsurprisingly not the only way to model data.

## 6.1 ORTHOGONAL POLYNOMIAL METHODS

Although the Measure Based Reconstruction is completely contained by working out two types of ergodic sum, the computational effort is quite high to work out the coefficients of the fitting. This high computational effort can lead to noise being added into the model just by the rounding error in the computer system. This exact problem was noted in the paper [5]. Their solution was to use monomials and perform a least square's fitting. These monomials could produce the orthogonal polynomials of the Measure based system, but the computation of the coefficients is significantly less due to the least square's method being a wide spread method. Their least square's problem is set out by

$$\Im = \sum_n \sum_{a=1}^{d_e} (X_a^{n+1} - F_a(\mathbf{X}^n))^2$$
$$\mathbf{F}(\mathbf{X}) \equiv \sum_k \mathbf{B}^k M_k(\mathbf{X})$$
(30)

and is computed respect to the $B^k$s where $M_k(X)$ represents all the monomials up to order $p$ that can be formed with the $d_e$ components $X_\alpha$ of the vectors $X$. The sum over $n$ runs over the whole data set.

This method has the advantage of not needing separate forms to calculate the multidimensional model. It also has the advantage of some well known minimisation algorithms.

## 6.2 Differential Equations.

As has been mentioned in Section 2 the underlying system of data set can be described using a differential equation, so it would seem logical that one would be able to model a data set with a differential equation. Briggs et al in their paper [12] present a method to perform such a modelling method. Basically the method is based on the idea that the invariant quantities of the data are not sufficient to fit the parameters of an ordinary differential equation. They the describe a method based upon the combination of a least squares fitting with an initial point problem. This is then used to adapt an error propagation method to compute the parameters for the ordinary differential equation.

## 6.3 Neural Networks and Stochastic Methods

As mention in the modelling section neural networks are powerful method for the modelling of data. As the area is quite new there are still a lot of conflicting reports on the subject and method to use. In general they are hard to set up, but once set up they do not expect much knowledge from the user. An example of this is shown in the paper [23]. These methods also show a similarity with that of the Measure Based Reconstruction method, in that after the training of the system the modelling can be taken from any point.

In terms of stochastic models, there has been great research into these methods as would be expected. Most papers on the subject present those method already discussed in this paper. Some suggest new methods for the process of parameter estimation, but generally they all agree that white noise is require for these methods.

# 7 Conclusion

In this project the subject has been presented. As modelling data is a large subject, the project has focussed on the main methods of modelling data and specifically the Measure Based Reconstruction Method. In the study of this method, the one dimensional model was tested in many numerical experiments, which proved that the method has a long prediction length at a very low error tolerance and is robust to noise. Unfortunately the two dimensional method had problems that meant that it could not be tested. However these problems were identified and explained with the specific part being worked out in terms of the lowest form. Also other methods were suggested as alternatives to the Measure based reconstruction and references where further information were given.

In general the application of a modelling method to real world data will not be as good as the application to the synthetic data in this project. This is because in most cases the order of the approximation is not always the optimal value for that variables. This might be due to many factors, some of which are not directly related to the model. In general the methods for determining the order of a model are sensitive to such factors as noise in the data or multiperiods. For some models there is no definite way of determining the optimal order and so it is more a matter of trial and error. Other factors that might affect the model in the real world are costs. These could be due to the time needed to run the method or the computation power required. These factors have not been investigated in this project as they are not related to the mathematics.

Further topics that have not been included in the project include such things as noise reduction of time series, neural networks and a greater investigation into other methods. References for these topics and all references used for this project are listed in the next section. These references are of about the same level as this papers and so should be easily understood.

[24] R. Hegger, H. Kantz, and T.Schreiber, **Practical implementation of nonlinear time series methods: The TISEAN package**, Chaos **9**:413, 1999

# A Project Log

In this section we shall be giving a log of the work done on the project, when it was done and if it doesn't make it into the project why it didn't. It will also give an idea of my thought processes, try to keep up!

## A.1 May

### A.1.1 May 1, 2002

I started to look for paper on the general methods of Modelling Chaotic Data such as Neural Networks and parameter fitting for differential equations. At this point it is only general research to see what area I would like to do, so it is basically reading papers and trying to understand the basics.

### A.1.2 May 12, 2002

Having found my supervisor (Dr Zaqueu Coelho) and discussed his ideas on what I should do as my dissertation topic, I decided at the end of the first full week of my initial research into the general area to write a progress report[7], with a list of the paper that I had read. In this report, which is included at the end of the project, I outlined my idea of a project description which combined the idea I had with those of Zaq.

### A.1.3 The rest of May, 2002

After my meeting in the beginning of May, I was waiting for two papers that Zaq was getting from Jason, and then photocopying them. While waiting for these papers, not much work could be done as I didn't exactly know what to research. So the basic idea of reading around the area seem to be the thing to do. In this time I also found a new MATLAB toolbox for neural networks that tied into a book on the subject. In this time I also tried to write note on the simpler method, such as local and global polynomial methods.

---

[7] This was mainly as Zaq told me that his computer was down, so it was just to confirm my e-mail.

## A.2 JUNE

### A.2.1 June 17, 2002

This was a good day as I received the papers from Zaq. By this time I had decided that I would photocopy the papers, and since Zaq suggested this idea as well it seemed like a good idea. The plan at this point was to read the paper, try and understand it and then go back and discuss what I was going to do with Zaq.

### A.2.2 June 19, 2002

After a general supervisor meeting with Jason, I am advised to read Chaos: the journal, as there are some recent papers in it that would be relevant to the project. So after the meeting I download the papers and read them. They are useful to the project.

Also I have a problem understanding the one dimensional method described in the paper [4]. Basically I think that the alpha coefficients will be needed to be calculated until negative infinity. Much working out on paper is done, but nothing seems to work.

### A.2.3 June 24, 2002

As part of the planning for the eventual programming of the method in one dimension I start to look into ways of generating data to test the methods against the claims in the paper [4]. I decided, as I am thinking of programming the method in MATLAB first, to use the TSTools toolbox to generate the Hénon data set. This toolbox is used as it was on the DAL[8] Disk that was provided by Jason in his course. Another thought at this time was to use the methods in the Tisean package, which is a command line based set of methods, and would be quicker for larger data sets. In other words if it is decided to write the methods in C, the Tisean package may be used.

---

[8] Data Analysis Laboratory

At this point there is still the problem with the alphas, but I am still trying things such as having the indices being modulo k. This seems to produce the fact that the problem alphas are equal to infinity, which is still not good. This is confirmed by Ian, who works it out separately to me and gets the same answer. Much paper is still being used trying to work it out.

## A.3 JULY

### A.3.1 July 5, 2002

Finally the problem with the indices has been resolved. Strangely it was solved on the Maths department stairs by Jason[9]. Basically I am reminded that linear functions do not have quadratic terms. By this I take it to mean that the alphas that are related to a order of x higher that a function of the order required at that point in the method will be zero. That is using $\alpha_i^j$, if j is greater than or equal to *i* then $\alpha$ has a value, otherwise it is zero. At this point the programming can start as the infinite loop has been removed.

### A.3.2 July 8, 2002

Initial programming of the algorithm is done in MATLAB, with another function where iterations of the functional approximation are included. At this point there is much cursing of MATLAB's array indexing as it goes from one instead of zero. For some reason the output looks strange, but I don't know why.

### A.3.3 July 9, 2002

Reprogramming of the algorithm with extra lines to get around the indexing problem. Basically the loops go from the values written in the paper, and then there is an additional number so that the arrays can be indexed. Still it seems to work better than before, but the results are not as expected.

---

[9] It's strange that it was on the stairs, not that Jason helped to solve it.

*A.3.4   July 10, 2002*

The full algorithm now works in MATLAB and seems to produce the same values as those that are computed by hand. When testing against the Hénon *x* variable with 50,000 points in the data set, the computation time is greater than 6 hours and so it is decided to write it in C. This is because C is a complied language while MATLAB's m-functions are interpreted. The programming into C is not done by using a converter but by me as I want to change certain items. This reprogramming takes until 5:00 but reduces the computational time to minutes instead of days.

*A.3.5   July 11, 2002*

The first test run of the program, and as always happens there is a problem. This is not surprising as I haven't programmed in C since last year, but that is just laziness. To start with it doesn't work under windows using my test Hénon *x*-variable. However after a bit of converting by hand the program compiles in UNIX (Tower) and then starts to work. Strange how these things happen like this. After much scratching of my head I decided to convert the windows program to use long doubles instead of doubles. This leads to more problems as it seems C isn't the same in Windows as it is in UNIX and for some reason that means sqrtl doesn't work in UNIX.

*A.3.6   July 12, 2002*

After the addition of the long doubles, the program now works under both operating systems, but a slight problem with the $N_k$s is discovered. This means accuracy is down. Iteration is added so that there are the number of iterations on each data point in the data file.

*A.3.7   July 14, 2002*

After staring at the code for two days I finally work out that squaring the sum of the variables is not the same as summing the squares of the variables. This is corrected in both programs and the accuracy suddenly shoots up.

However the graph of the Hénon attractor still shows a "box" that seems to be there after the second iteration.

### A.3.8   July 16, 2002

After not getting the results that I want, I try changing the algorithm so that the whole algorithm is iterated. This produces a cone shaped set of Hénon attractors and so is changed back.

### A.3.9   July 17, 2002

Saw Jason after the Board of Studies meeting and discuss the method of prediction to be used, as I had two failed attempts. The notes used to explain the method to be utilise are shown below:

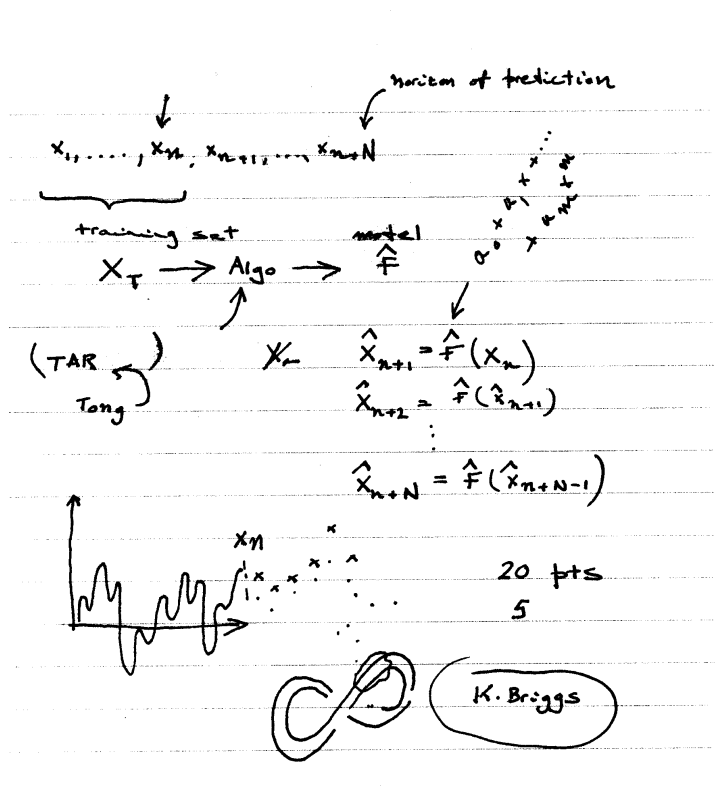

After this discussion the iterations are performed on the last data point of the data file, and then on each predicted point until the required amount are produced. Also Jason says that I should look at the NTL toolbox for C++ on Keith Briggs' website. This is done on my return to the MSc Room along with the start of the log[10].

### A.3.10 July 18, 2002

The changes work in both C and MATLAB versions of the program. As there are now fewer calculations, execution times fall. Also at this point I start to think about the experimentation that will be performed, and write a data generation program in C. I decided to write it as I can not find a program to create quadratic logistic data, or exponential quadratic logistic data. It is written in C as the data sets to be created will be 30,000+ points and I think MATLAB takes too long. Also as it is my own data program I can add more data set too it.

Also I take a look at the manual for the NTL toolbox. It seems to produce a greater degree of accuracy, which at this time I do not need as it is already allocating over 1 MB of memory to the variables. Also on Jason's advice I start to look at converting to Mex functions in MATLAB so that I can use C and have the variable imported into MATLAB directly. After reading many files I decided to stick with what is working.

### A.3.11 July 20, 2002

Added the one dimensional and two dimensional Hénon data generators into the data generating program. Also I start to look at the two dimensional method described in [4]. As the method looks to be of a greater difficulty I start by programming the moment generating functions.

### A.3.12 July 21, 2002

Worked out loop order for the two dimensional method and found a problem in the method. Also worked out that A and B from the two dimensional method are equivalent to b in the one dimensional method, and a and g form the two dimensional method are equivalent to a in the one dimensional method. So the two dimensional model has a similar structure to the one dimensional method. This means that the programming order has already be decided.

---

[10] So up to now it has been from memory, that's why some days are not included. I'm sure I was working on them, but I just can remember on what!

### A.3.13  July 22, 2002

Added exponential quadratic data generation to the infamous data generation program. Also finished programming the two dimensional method in C. This is not tested as the problem with the method, which is that an index is missing in the generation of the $a$'s, is still unresolved. Coding mirrors the structure of the paper, which may not be the most efficient way of doing things. However when something works there is no need to fix it.

### A.3.14  July 23, 2002

Reprogrammed the method to generated Hénon data so that it is now computed in the two dimensional subroutine rather than trying to do it in the one dimensional subroutine. Also fixed some of the errors in the two dimensional program, with a dummy index in the missing spot. Also corrected some errors in the data generation program. Also I find an initial point to generate quadratic data, without periodicity within 60,000 points, which is 0.123456789.

### A.3.15  July 25, 2002

As the missing index problem is quite important, I talk to Jason about it, and the fact that I have tried all the logical indexes and none seem to produce accurate prediction or they produce accurate predictions but they produce error in the program. Back in the one dimensional method I try to prove that the graph in the paper [4] is correct. This turns out to be true but by not using their method to work out the prediction length.

Also I try to work out by hand why the two dimensional method doesn't work. Also this is a hope to find out which index it might be. I go back to working it out from the moments of the data, but it still produces a negative $N_{2,0}$ so it cannot be square rooted with a real result. However it will work if $\eta^{(i,j)} = 0$, but that only happens if the missing index is i. This will produce a quite bad approximation and so is rejected as a solution.

### A.3.16  July 26, 2002

Still there is a slight problem in the programming of the two dimensional method, but I still can't find it. After more than three hours of debugging in one day I decided that it is the missing index which causes all of the problems. After talking again with Jason, it is decided to e-mail one of the authors of the paper to see if there is a canned routine. As a quite extensive search of many websites only produces a link for Massimiliano Giona. So I e-mail him to see if the missing index problem can be solved, and if there is a canned routine.

### A.3.17  July 30, 2002

I prepare for a meeting with Zaq tomorrow by printing out the work that has currently been done. Also I run the programs with quadratic and exponential quadratic data to produce more graphs for the meeting. Most do not work or do not show what is required so are rejected.

### A.3.18  July 31, 2002

I have a meeting with Zaq and discuss how the project is going and what has been done. Also we decided that the Hénon data set that was used to check the one dimensional method should be ignore as there is no direct relationship between consecutive points. He suggest that the missing index would either be h or k, but I already know that it will not work with those. Also I have another go at debugging the two dimensional method with little success except to find that the value that breaks the program and algorithm is a factor of ten out depending on how it is calculated. Overall it is still negative.

## A.4 AUGUST

### A.4.1 August 1, 2002

After noting a great error in my working for the two dimensional problem, I rewrote the working fixing the error. Also Jason discussed the project and suggested some papers to read to round the project out. Also he told me that the problem that I have with the second paper [5] is because the missing coefficients are supposed to be the item that I am calculating. I then download the paper that he has recommended and look again at the second method.

### A.4.2 August 2, 2002

Today I battle with the University's interlibrary loan system, eventually going back to basics and filling out a paper form to get the book [1].

### A.4.3 August 3, 2002

I finally get round to reading the papers that Jason suggested that I should read, and try and work out how to program the listing that is given in [10].

### A.4.4 August 5, 2002

Today I add the prediction time calculation as it is given in the paper [4]. As far as I can see it is just getting the average prediction time of all the points in the data set.. This makes the programs go slightly slower again but it does eventually produce some results. However this does not seem to produce a graph similar to that given in the paper. Also the graphs that were hinting at by Jason in the meeting on the 1st are produced for the quadratic data.

### A.4.5 August 6, 2002

Today I tried to program both the least squares method and the method given in the paper [10]. Neither of these seem to work or produce results that seem to be of any use. So I go searching on the web for some addition information on the method from the paper, but this just seems to be wasting time. The graphs for the exponential quadratic are produced.

*A.4.6 August 12, 2002*

Tried to program additive noise to the data generation program but fail due to there not being a Gaussian random number generator in either C or C++. Other methods to add noise to data are investigated including importing into MATLAB, adding noise and the exporting the file.

*A.4.7 August 14, 2002*

Using the TISEAN package addnoise program the noise is added to the data sets. After finally giving up on the two dimensional MB model and the least squares method, all experiments are performed using the one dimensional model.

*A.4.8 August in general.*

Throughout the whole of August the project is written up. This takes many drafts and many different directions are investigated, but few are chosen. Also the periodogram of the data sets are drawn to prove that the data's properties are being carried through the models. The noise values are investigated, and the regulations of the format of the document are changed. In general each day is filled with large amounts of work and most days do not finish till late at night, due to the writing up.

A.5 SEPTEMBER

*A.5.1 September 2, 2002*

The project is handed in!

# B PROGRAM LISTINGS

The program listings for the programs used in the project now follow. In general they are written in the C programming language, but a few exceptions are written in MATLAB m-files.